\documentclass[b5paper,10pt]{article}
\usepackage{amssymb}
\usepackage{graphicx}
\usepackage{graphics}
\def\yen{\hbox{iftdir\yoko\fi

\setbox0=\hbox{Y}Y\kern-.97\wd0\vbox{\hrule height.lex width.98\wd0
\kern.33ex\hrule height.lex width.98\wd0\kern.45ex}}}

\def\yen{{\setbox0=\hbox{Y}Y\kern-.97\wd0\vbox{hrule height.lex width.98%
\wd0\kern.33ex\hrule height.lex width.98\wd0\kern.45ex}}}

\def\x{\times }
\def\f{\newline }
\begin{document}

\baselineskip=12pt

\title{
Nonribbon 2-links \\
all of whose components are trivial knots and \\
some of whose band-sums are nonribbon knots 
}     
\author{ Eiji Ogasa 
\footnote{
This research was partially supported 
by Research Fellowships
of the Promotion of Science for Young Scientists.
\newline
{\bf Keyword:}
2-links, 2-knots, nonribbon 2-links, nonribbon 2-knots
\newline
{\it 1991 Mathematics Subject Classification.} 
Primary 57M25, 57Q45 
\newline
This paper is published in \newline
Journal of knot theory and its ramificatioms 10, 2001, 913-922.\newline
This manuscript is not the published version.
} \\
 ogasa@hep-th.phys.s.u-tokyo.ac.jp\\
High Energy Physics Theory Group\\ 
Department of Physics\\ 
 University of Tokyo\\ 
Hongo, Tokyo 113, JAPAN\\
}
\date{}
\maketitle

\noindent
{\bf Abstract.}
There is a nonribbon 2-link all of whose components are trivial 2-knots 
and one of whose band-sums is a nonribbon 2-knot.

\section{ Main result }

We work in the smooth category. 

An {\it $m$-component 2-(dimensional) link} 
is a closed oriented 2-submanifold 
$L=(K_1,...,K_m)$ $\subset S^4$ 
such that $K_i$ is diffeomorphic to $S^2$. 
If $m=1$, $L$ is called a {\it 2-knot} 
We say that 2-links $L_1$ and $L_2$ are {\it equivalent} 
if there exists an orientation preserving diffeomorphism 
$f:$ $S^4$ $\rightarrow$ $S^4$ 
such that $f(L_1)$=$L_2$  and 
that $f | _{L_1}:$ $L_1$ $\rightarrow$ $L_2$ is 
an order and orientation preserving diffeomorphism.    
Take a 3-ball $B^3$ in $S^4$. Then $\partial B^3$ is a 2-knot. 
We say that a 2-knot $K$ is a {\it trivial} knot if 
$K$ is equivalent to the 2-knot $\partial B^3$.

A 2-link $L= ( K_1,...,K_m ) $ is called a {\it ribbon } 2-link 
if $L$ satisfies the following properties. 
(See e.g. 
\cite{C}.)
\f
(1) There is a self-transverse immersion 
$f:D^3_1\amalg...\amalg D^3_m$
$\rightarrow S^4$ 
such that $f(\partial D^3_i)=K_i$.  
\f
(2) The singular point set $C$  $(\subset S^4$) 
of $f$ consists of double points. 
$C$ is a disjoint union of 2-discs $D^2_i (i=1,...,k)$. 
\f
(3) Put $f^{-1}(D^2_j)=D^2_{jB}\amalg D^2_{jS}$. 
The 2-disc 
$D^2_{jS}$ is trivially embedded in the interior Int $D^3_\alpha$ 
of a 3-disc component $D^3_\alpha$.  
The circle $\partial D^2_{jB}$ is trivially embedded in 
the boundary $\partial D^3_\beta$ of a 3-disc component $D^3_\beta$. 
The 2-disc $D^2_{jB}$ is trivially embedded 
in the 3-disc component $D^3_\beta$. 
(Note that there are two cases, $\alpha=\beta$ and $\alpha\neq\beta$.)

There are nonribbon 2-knots. (See e.g. 
\cite{C}.)  
It is trivial that, if a component of a 2-link is a nonribbon 2-knot, 
the 2-link is a nonribbon 2-link. It is natural to ask:  

\vskip3mm
\noindent
{\bf Question}
Is there a nonribbon 2-link all of whose components are ribbon knots?
In particular, 
is there a nonribbon 2-link all of whose components are trivial knots?   

\vskip3mm

We give an affirmative answer to this question. 

\vskip3mm
\noindent
{\bf Theorem 1.1 }
{\it 
There is a nonribbon 2-link $L=(K_1, K_2)$ 
such that $K_i$ is a trivial 2-knot 
($i=1,2$).  }

\vskip3mm
\noindent
{\bf Note.} 
The announcement of Theorem 1.1 is in 
\cite{O1}.

\section{ Band-sums}

Let $L=(K_1, K_2)$ be a 2-link. 
A 2-knot $K_0$  is called 
a {\it band-sum } of the components $K_1$ and $K_2$ of
the 2-link $L$ along a {\it band} $h$ 
if we have: 
\f
(1) There is a 3-dimensional 1-handle $h$, 
which is attached to $L$, embedded in $S^4$. 
\f
(2) 
There are a point $p_1\in K_1$ and 
a point $p_2\in K_2$. 
We attach $h$ to $K_1\amalg K_2$ along $p_1\amalg p_2$. 
 $h\cap (K_1\cup K_2)$ is the attach part of $h$.  
Then we obtain a 2-knot from $K_1$ and $K_2$ 
by this surgery. The 2-knot is $K_0$.

\section{ 
A sufficient condition of Theorem 1.1  
}

In \S4 and \S5 we prove: 

\vskip3mm
\noindent
{\bf Proposition 3.1 }{\it 
There is a 2-link $L=(K_1,K_2)$ such that 
\f
(1) $K_i$ is a trivial 2-knot ($i=1,2$), and 
\f
(2) a band-sum $K_3$ of 
the components $K_1,K_2$ of the 2-link 
$L$ is a nonribbon 2-knot. }

\vskip3mm
\noindent
{\bf 
 Claim 3.2  }{\it 
Proposition 3.1 implies Theorem 1.1. }

\vskip3mm
\noindent
{\bf Proof of Claim 3.2. }  
By the definition of ribbon links, we have the following fact: 
If $L=(K_1, K_2)$ is a ribbon 2-link, then any band sum 
of $L=(K_1,K_2)$  is a ribbon 2-knot. 
The contrapositive proposition of this fact implies Claim 3.2.

\section{  $Q(K)$ }

Let $K$ be a 2-knot $\subset S^4$. We define a 2-knot $Q(K)$ for $K$. 
The 2-knot $Q(K)$ plays an important role in our proof 
as we state in the last paragraph of this section. 

Let $K\x D^2$ be a tubular neighborhood of $K$ in $S^4$, 
where $D^2$ is a disc. 
In Int $D^2$, take a compact oriented 1-dimensional submanifold $[-1,1]$. 
Take  $K\x[-1,1]$ $\subset K\x$ Int$D^2$. 
We give an orientation to  $K\x[-1,1]$. 
Let $D(K)$ be the 2-component 2-link 
$(K\x\{-1\}, K\x\{1\})$. 
Then $K\x[-1,1]$ is a Seifert hypersurface of the 2-link $D(K)$, 
where we give an orientation to $D(K)$  so that 
the orientation of $D(K)$ is compatible with that of  $K\x[-1,1]$.

In order to prove our main theorem,  
we construct some 2-knots, 2-links, and some subsets in $S^4$ from $D(K)$. 
For this purpose, we prepare the following $B^4$ and $F_\theta$. 
\f Let $B^4$ be a 4-ball $\subset S^4$. 
Put $B^4=$

$\{( x,y,z,w )\vert$
$0\leq x \leq1, 0\leq y \leq1, $
$z=r\cdot cos\theta, w=r\cdot sin\theta, $
$0\leq r \leq1, 0\leq\theta<2\pi \}$. 
\f Let $F _0=$
$\{( x,y,z,0 )\vert$
$0\leq x \leq1, 0\leq y \leq1, 0\leq z \leq1 \}\subset B^4$. 
\f Let 
$A=\{( x,y,0,0 )\vert$
$0\leq x \leq1, 0\leq y \leq1 \}\subset F_0\subset B^4$. 
\f We regard $B^4$ as the result of rotating 
$F_0$ around the axis $A$. 
For each $\theta$, we put 

$F_\theta=$
$\{(x,y,r\cdot cos\theta, r\cdot sin\theta)\vert$
$0\leq x \leq1, 0\leq y \leq1, 0\leq r \leq1,\theta$: 
fix$\}$.

We suppose that $B^4\cap D(K)$ satisfies the condition that, 
 for each $\theta$,  $F_\theta \cap D(K)$ is drawn as in Figure 4.1.

\newpage
\unitlength 0.1in
\begin{picture}(53.80,48.20)(2.00,-47.80)
%
\special{pn 8}%
\special{pa 1000 790}%
\special{pa 1000 4780}%
\special{fp}%
%
\special{pn 8}%
\special{pa 2010 30}%
\special{pa 1000 810}%
\special{fp}%
%
\special{pn 8}%
\special{pa 2010 3990}%
\special{pa 1000 4770}%
\special{fp}%
%
\special{pn 8}%
\special{pa 1600 1080}%
\special{pa 3700 1080}%
\special{fp}%
%
\special{pn 8}%
\special{pa 1610 800}%
\special{pa 4010 800}%
\special{fp}%
\put(7.9000,-47.1000){\makebox(0,0)[lb]{y}}%
\put(20.9000,-1.3000){\makebox(0,0)[lb]{x}}%
\put(19.0000,-3.6000){\makebox(0,0)[lb]{1}}%
\put(10.7000,-45.9000){\makebox(0,0)[lb]{1}}%
\put(54.0000,-42.1000){\makebox(0,0)[lb]{1}}%
\put(20.1000,-41.5000){\makebox(0,0)[lb]{o}}%
%
\special{pn 8}%
\special{pa 3710 1090}%
\special{pa 3710 1280}%
\special{fp}%
%
\special{pn 8}%
\special{pa 3710 1540}%
\special{pa 3710 2390}%
\special{fp}%
%
\special{pn 8}%
\special{pa 3710 4360}%
\special{pa 3710 2900}%
\special{fp}%
%
\special{pn 8}%
\special{pa 4040 3780}%
\special{pa 4040 4360}%
\special{fp}%
%
\special{pn 8}%
\special{pa 4050 3530}%
\special{pa 4050 2910}%
\special{fp}%
%
\special{pn 8}%
\special{pa 4050 2380}%
\special{pa 4040 800}%
\special{fp}%
%
\special{pn 8}%
\special{pa 3980 810}%
\special{pa 4040 800}%
\special{fp}%
%
\special{pn 8}%
\special{pa 380 650}%
\special{pa 390 650}%
\special{fp}%
%
\special{pn 8}%
\special{pa 200 1200}%
\special{pa 200 1210}%
\special{fp}%
%
\special{pn 4}%
\special{pa 1000 200}%
\special{pa 960 210}%
\special{dt 0.027}%
\special{pa 960 210}%
\special{pa 961 210}%
\special{dt 0.027}%
\put(55.8000,-39.2000){\makebox(0,0)[lb]{r}}%
\put(25.3000,-47.6000){\makebox(0,0)[lb]{Figure 4.1}}%
%
\special{pn 8}%
\special{pa 3700 1250}%
\special{pa 3710 1600}%
\special{fp}%
%
\special{pn 8}%
\special{pa 3730 2360}%
\special{pa 3710 3000}%
\special{fp}%
%
\special{pn 8}%
\special{pa 4050 2290}%
\special{pa 4050 2960}%
\special{fp}%
%
\special{pn 8}%
\special{pa 4050 3530}%
\special{pa 4040 3850}%
\special{fp}%
\put(34.0000,-22.0000){\makebox(0,0)[lb]{E}}%
\put(42.6000,-22.2000){\makebox(0,0)[lb]{E}}%
\put(35.0000,-22.9000){\makebox(0,0)[lb]{1}}%
\put(43.7000,-23.0000){\makebox(0,0)[lb]{2}}%
%
\special{pn 8}%
\special{pa 2020 30}%
\special{pa 2020 690}%
\special{fp}%
%
\special{pn 8}%
\special{pa 2020 900}%
\special{pa 2020 1030}%
\special{fp}%
%
\special{pn 8}%
\special{pa 2020 1220}%
\special{pa 2020 3990}%
\special{fp}%
%
\special{pn 8}%
\special{pa 2020 3970}%
\special{pa 3580 3970}%
\special{fp}%
%
\special{pn 8}%
\special{pa 3780 3970}%
\special{pa 3960 3970}%
\special{fp}%
%
\special{pn 8}%
\special{pa 4160 3970}%
\special{pa 5560 3970}%
\special{fp}%
\end{picture}%

\vskip10mm
\noindent
{\bf Note.} In Figure 4.1, we suppose the following hold: 
The intersection $B^4\cap D(K)$ is 
a disjoint union of two 2-discs. 
Call them $D^2_1$ and $D^2_2$. 
The intersection $F_\theta \cap D(K)$ is two arcs. 
Call them $E_1$ and $E_2$. 
The boundary $\partial E_i$ is a set of two points $a_i\amalg b_i$, 
where $a_i$ is in $A$ and $b_i$ is in $F_0-A$. 
The 2-disc $D_i$ is the result of rotating $E_i$ around the axis $A$. 
The result of rotating $b_i$ is $\partial D^2_i$.
Since $a_i$ is in the axis $A$, 
the result of rotating $a_i$ is the point $a_i$ itself. 
The point $b_i$ is in the boundary of $D^2_i$.  
The point $a_i$ is in the interior of $D^2_i$.  
\vskip3mm

Let $Q(K)$ be a band-sum of the components $K\x\{-1\}$ and $K\x\{1\}$ 
of the 2-link $D(K)$ with the following properties.  
\f
(1)
 The band $h$ is in $K\x$ Int$D^2$. 
\f
(2)    $\{h-$(the attach part of $h$)$\}$$\cap (K\x[-1,1])=\phi$. 
\f
(3)    $\overline{Q(K)-B^4}=$ $\overline{D(K)-B^4}$. 
\f
(4)
 $B^4\cap$  $(D(K)\cup h)$
 ( =$B^4\cap$ $(K\x\{-1\}\cup h\cup K\x\{1\})$ ) 
satisfies the following conditions. 
(We summarize the conditions in Table 1.)
\f
For $\pi\leq\theta<2\pi$ and $\theta=0$, 
$F_\theta\cap( D(K)\cup h)$ is drawn as in Figure 4.2. 
\f
For $0<\theta<\pi$, 
$F_\theta\cap(D(K)\cup h)$ is drawn as in Figure 4.1.
\f
(5)
 $B^4\cap h$ satisfies the following conditions.  
\f
For $\pi\leq\theta<2\pi$ and $\theta=0$, 
$F_\theta\cap h$ is drawn as in Figure 4.3. 
\f
For $0<\theta<\pi$, 
$F_\theta\cap h$  is empty. 
\f
$B^4\cap$ (the attach part of $h$)  is as follows.  
\f
For $\pi\leq\theta<2\pi$ and $\theta=0$, 
$F_\theta\cap$ (the attach part of $h$) is drawn as in Figure 4.4. 
\f
For $0<\theta<\pi$, 
$F_\theta\cap$ (the attach part of $h$) is empty. 
\f
(6)
$B^4\cap Q(K)$ satisfies the following conditions. 
(We summarize the conditions in Table 1.)
\f
For $\pi<\theta<2\pi$, $F_\theta\cap Q(K)$ is drawn as in Figure 4.5. 
\f
For $\theta=0,\pi$, $F_\theta\cap Q(K)$ is drawn as in Figure 4.2.  
\f
For $0<\theta<\pi$, $F_\theta\cap Q(K)$ is drawn as in Figure 4.1.

\newpage
\input 4.2.tex

\newpage
\input 4.3.tex

\newpage
\unitlength 0.1in
\begin{picture}(53.70,48.20)(2.00,-47.80)
%
\special{pn 8}%
\special{pa 2010 10}%
\special{pa 2010 4000}%
\special{fp}%
%
\special{pn 8}%
\special{pa 1000 790}%
\special{pa 1000 4780}%
\special{fp}%
%
\special{pn 8}%
\special{pa 2010 30}%
\special{pa 1000 810}%
\special{fp}%
%
\special{pn 8}%
\special{pa 2010 3990}%
\special{pa 1000 4770}%
\special{fp}%
%
\special{pn 8}%
\special{pa 2010 3990}%
\special{pa 5530 3990}%
\special{fp}%
\put(7.9000,-47.1000){\makebox(0,0)[lb]{y}}%
\put(20.9000,-1.3000){\makebox(0,0)[lb]{x}}%
\put(19.0000,-3.6000){\makebox(0,0)[lb]{1}}%
\put(10.7000,-45.9000){\makebox(0,0)[lb]{1}}%
\put(54.0000,-42.1000){\makebox(0,0)[lb]{1}}%
\put(20.1000,-41.5000){\makebox(0,0)[lb]{o}}%
%
\special{pn 8}%
\special{pa 380 650}%
\special{pa 390 650}%
\special{fp}%
%
\special{pn 8}%
\special{pa 200 1200}%
\special{pa 200 1210}%
\special{fp}%
%
\special{pn 4}%
\special{pa 1000 200}%
\special{pa 960 210}%
\special{dt 0.027}%
\special{pa 960 210}%
\special{pa 961 210}%
\special{dt 0.027}%
%
\special{pn 4}%
\special{pa 4050 3540}%
\special{pa 4050 3800}%
\special{dt 0.027}%
\special{pa 4050 3800}%
\special{pa 4050 3799}%
\special{dt 0.027}%
%
\special{pn 8}%
\special{pa 3720 1520}%
\special{pa 3710 1290}%
\special{dt 0.045}%
\special{pa 3710 1290}%
\special{pa 3710 1291}%
\special{dt 0.045}%
%
\special{pn 8}%
\special{pa 3710 1290}%
\special{pa 3730 1540}%
\special{fp}%
%
\special{pn 8}%
\special{pa 4060 3800}%
\special{pa 4040 3520}%
\special{fp}%
\put(55.7000,-38.9000){\makebox(0,0)[lb]{r}}%
\put(22.3000,-47.2000){\makebox(0,0)[lb]{Figure 4.4}}%
\end{picture}%

\newpage
\unitlength 0.1in
\begin{picture}(54.40,48.20)(2.00,-47.80)
%
\special{pn 8}%
\special{pa 1000 790}%
\special{pa 1000 4780}%
\special{fp}%
%
\special{pn 8}%
\special{pa 2010 30}%
\special{pa 1000 810}%
\special{fp}%
%
\special{pn 8}%
\special{pa 2010 3990}%
\special{pa 1000 4770}%
\special{fp}%
%
\special{pn 8}%
\special{pa 1600 1080}%
\special{pa 3700 1080}%
\special{fp}%
%
\special{pn 8}%
\special{pa 1610 800}%
\special{pa 4010 800}%
\special{fp}%
\put(7.9000,-47.1000){\makebox(0,0)[lb]{y}}%
\put(20.9000,-1.3000){\makebox(0,0)[lb]{x}}%
\put(19.0000,-3.6000){\makebox(0,0)[lb]{1}}%
\put(10.7000,-45.9000){\makebox(0,0)[lb]{1}}%
\put(54.0000,-42.1000){\makebox(0,0)[lb]{1}}%
\put(20.1000,-41.5000){\makebox(0,0)[lb]{o}}%
%
\special{pn 8}%
\special{pa 3710 1280}%
\special{pa 2690 1280}%
\special{fp}%
%
\special{pn 8}%
\special{pa 3710 1530}%
\special{pa 2910 1530}%
\special{fp}%
%
\special{pn 8}%
\special{pa 2700 1290}%
\special{pa 2700 2140}%
\special{fp}%
%
\special{pn 8}%
\special{pa 2920 1530}%
\special{pa 2920 1890}%
\special{fp}%
%
\special{pn 8}%
\special{pa 5110 1910}%
\special{pa 5110 2750}%
\special{fp}%
%
\special{pn 8}%
\special{pa 4820 2130}%
\special{pa 4820 2490}%
\special{fp}%
%
\special{pn 8}%
\special{pa 4820 2510}%
\special{pa 2700 2500}%
\special{fp}%
%
\special{pn 8}%
\special{pa 2710 2500}%
\special{pa 2710 3260}%
\special{fp}%
%
\special{pn 8}%
\special{pa 5110 2760}%
\special{pa 2920 2760}%
\special{fp}%
%
\special{pn 8}%
\special{pa 2940 2760}%
\special{pa 2940 3050}%
\special{fp}%
%
\special{pn 8}%
\special{pa 5110 3050}%
\special{pa 5110 3750}%
\special{fp}%
%
\special{pn 8}%
\special{pa 4820 3290}%
\special{pa 4820 3530}%
\special{fp}%
%
\special{pn 8}%
\special{pa 4820 3510}%
\special{pa 4020 3530}%
\special{fp}%
%
\special{pn 8}%
\special{pa 5110 3780}%
\special{pa 4020 3790}%
\special{fp}%
%
\special{pn 8}%
\special{pa 3710 1090}%
\special{pa 3710 1280}%
\special{fp}%
%
\special{pn 8}%
\special{pa 3710 1540}%
\special{pa 3710 2390}%
\special{fp}%
%
\special{pn 8}%
\special{pa 3710 4360}%
\special{pa 3710 2900}%
\special{fp}%
%
\special{pn 8}%
\special{pa 3710 2570}%
\special{pa 3710 2680}%
\special{fp}%
%
\special{pn 8}%
\special{pa 4040 3780}%
\special{pa 4040 4360}%
\special{fp}%
%
\special{pn 8}%
\special{pa 4050 3530}%
\special{pa 4050 2910}%
\special{fp}%
%
\special{pn 8}%
\special{pa 4060 2680}%
\special{pa 4050 2580}%
\special{fp}%
%
\special{pn 8}%
\special{pa 4050 2380}%
\special{pa 4040 800}%
\special{fp}%
%
\special{pn 8}%
\special{pa 3980 810}%
\special{pa 4040 800}%
\special{fp}%
%
\special{pn 8}%
\special{pa 2920 1880}%
\special{pa 3570 1880}%
\special{fp}%
%
\special{pn 8}%
\special{pa 3800 1880}%
\special{pa 3970 1880}%
\special{fp}%
%
\special{pn 8}%
\special{pa 4190 1890}%
\special{pa 5120 1910}%
\special{fp}%
%
\special{pn 8}%
\special{pa 2700 2140}%
\special{pa 3570 2140}%
\special{fp}%
%
\special{pn 8}%
\special{pa 3800 2130}%
\special{pa 3960 2120}%
\special{fp}%
%
\special{pn 8}%
\special{pa 4190 2120}%
\special{pa 4830 2140}%
\special{fp}%
%
\special{pn 8}%
\special{pa 2930 3040}%
\special{pa 3600 3040}%
\special{fp}%
%
\special{pn 8}%
\special{pa 3800 3040}%
\special{pa 3980 3040}%
\special{fp}%
%
\special{pn 8}%
\special{pa 4170 3030}%
\special{pa 5110 3060}%
\special{fp}%
%
\special{pn 8}%
\special{pa 380 650}%
\special{pa 390 650}%
\special{fp}%
%
\special{pn 8}%
\special{pa 200 1200}%
\special{pa 200 1210}%
\special{fp}%
%
\special{pn 8}%
\special{pa 2710 3260}%
\special{pa 3600 3270}%
\special{fp}%
%
\special{pn 8}%
\special{pa 3800 3290}%
\special{pa 3970 3300}%
\special{fp}%
%
\special{pn 8}%
\special{pa 4170 3280}%
\special{pa 4820 3310}%
\special{fp}%
\put(41.8000,-7.9000){\makebox(0,0)[lb]{l}}%
\put(42.5000,-8.5000){\makebox(0,0)[lb]{1}}%
\put(52.0000,-35.9000){\makebox(0,0)[lb]{l}}%
\put(52.6000,-36.7000){\makebox(0,0)[lb]{2}}%
\put(22.2000,-47.7000){\makebox(0,0)[lb]{Figure 4.5}}%
\put(55.8000,-39.2000){\makebox(0,0)[lb]{r}}%
%
\special{pn 8}%
\special{pa 2020 50}%
\special{pa 2020 660}%
\special{fp}%
%
\special{pn 8}%
\special{pa 2020 870}%
\special{pa 2020 1030}%
\special{fp}%
%
\special{pn 8}%
\special{pa 2020 1190}%
\special{pa 2020 3990}%
\special{fp}%
%
\special{pn 8}%
\special{pa 2020 3990}%
\special{pa 3600 3990}%
\special{fp}%
%
\special{pn 8}%
\special{pa 3790 3990}%
\special{pa 3980 3990}%
\special{fp}%
%
\special{pn 8}%
\special{pa 4140 3990}%
\special{pa 5640 3990}%
\special{fp}%
\end{picture}%

\newpage

\begin{tabular}{|l||l|l|l|l|}
   &  {$\theta=0$} & {$0<\theta<\pi$} & {$\theta=\pi$} & {$\pi<\theta<2\pi$}\\ 
 $B^4\cap$$D(K)$      
&Figure 4.1 & Figure 4.1 & Figure 4.1&Figure 4.1   \\
$B^4\cap (D(K)\cup h)$  
&Figure 4.2 & Figure 4.1 & Figure 4.2 &Figure 4.2   \\
$B^4\cap Q(K)$ 
&Figure 4.2 & Figure 4.1 & Figure 4.2 &Figure 4.5   \\
$B^4\cap L$  
&Figure 4.5 & Figure 4.5 & Figure 4.5 &Figure 4.5   \\
 $B^4\cap$  $(L\cup h')$ 
&Figure 4.2 & Figure 4.2 & Figure 4.2 &Figure 4.5   \\
$B^4\cap K_3$ 
&Figure 4.2 & Figure 4.1 & Figure 4.2 &Figure 4.5   \\
\end{tabular}

\vskip1cm
\hskip5cm Table 1

\noindent
{\bf Note.} The following hold: 
\f
(I) $h\cup (K\x[-1,1])$ is a Seifert hypersurface of the 2-knot $Q(K)$.    
\f
(II) Let $A$ be a Seifert matrix of $Q(K)$ associated with  
the Seifert hypersurface $h\cup (K\x[-1,1])$. 
Since $h\cup (K\x[-1,1])$ is diffeomorphic to $\overline{(S^1\x S^2)-B^3}$,  
$A$ is a $(1\x1)$-matrix. 
By the construction of $Q(K)$, $A=(2)$ or $A$=(-1) holds. 
Recall that whether $A$ is a $1\times1$-matrix (2) or (-1) depends on which orientation we give $Q(K)$. 
Recall that the orientation of $Q(K)$ is determined by that of $D(K)$.
\f
(III) $2t-1$ or $2-t$ represents for the Alexander polynomial of the 2-knot $Q(K)$.   
(See \S F,G,H of \S 7 
of 
\cite{Ro} 
for Seifert matrices of 2-knots and the Alexander polynomial of 2-knots. )
Hnece $Q(K)$ is a nontrivial 2-knot.  
\vskip3mm

In \S5 we prove:

\vskip3mm
\noindent
{\bf Lemma 4.1 }{\it 
 Let $K$ be a 2-knot. There is a 2-link $L=(K_1, K_2)$ such that 
\f
(1) $K_i$ is a trivial 2-knot ($i=1,2$), and   
\f
(2) $Q(K)$ is a band-sum $K_3$ of the components $K_1$, $K_2$ of 
the 2-link $L$.}
\vskip3mm

The above $Q(K)$ is `a 2-knot $D(J,\gamma)$ whose  $\gamma$ 
is sufficiently complicated' 
in \S4 of 
\cite{C}. 
Corollary 4.3 in \S4 of 
\cite{C} 
or  Example after Corollary 4.3 in \S4 of 
\cite{C} 
says that, for a 2-knot $K$, the above $Q(K)$ is a nonribbon 2-knot.  
Hence Lemma 4.1 implies Proposition 3.1.

\section{  Proof of Lemma 4.1 }

Let $L=(K_1, K_2)$ be a 2-link with the following conditions.  
\f
(1) $(S^4-B^4)\cap L$=$(S^4-B^4)\cap D(K)$. 
\f
(2)
 $B^4\cap L$ satisfies the condition that, 
for each $\theta$, $F_\theta\cap L$ is drawn as in Figure 4.5. 
(We summarize the conditions in Table 1.)

\vskip3mm
\noindent
{\bf Note.} 
In Figure 4.5, the following hold: 
The two arcs are called $l_1$ and $l_2$. 
 $l_i$ is a trivial arc.  
 $l_1\cap A\neq\phi$.  
 $l_2\cap A=\phi$.  
 $K_i$ is made from $l_i$ by the rotation.
\vskip3mm

By the construction of $L=(K_1, K_2)$, $K_1$ satisfies the conditions: 
\f
(1) $K_1\subset B^4$. 
\f
(2) For each $\theta$, $F_\theta\cap K_1$ is drawn as in Figure 5.1.   

We prove $K_1$ is a trivial knot. 
Because: Since $l_1$ is a trivial arc, $K_1$ is a spun knot of a trivial 1-knot. 
See 
\cite{Z}  
for spun knots.

By the construction of $L=(K_1, K_2)$, $K_2$ satisfies the following conditions. \f
(1) $(S^4-B^4)\cap K_2$=$(S^4-B^4)\cap D(K)$. 
\f
(2) For each $\theta$, $F_\theta\cap K_2$ is drawn as in Figure 5.2.

We prove: $K_2$ is a trivial knot. 
Because: 
Let $P$ be a subset $(K\x[-1,1])-B^4\}$. Then $P$ is diffeomorphic to a 3-ball. 
Hence  $\partial P$ is a trivial 2-knot. 
Since $l_2$ is a trivial arc, $K_2$ is equivalent to $\partial P$. 
Hence $K_2$ is a trivial 2-knot.

Let  $K_3$ be a band-sum of the components $K_1$ and $K_2$ of the 2-link $L$ 
with the following conditions. 
\f
(1) The band $h'$ is in $B^4$. 
\f
(2) $B^4\cap$$(L\cup h')$=$B^4\cap(K_1\cup h' \cup K_2)$ 
satisfies the following conditions. (We summarize the conditions in Table 1.)
\f
For $\pi<\theta<2\pi$, $F_\theta\cap(L\cup h')$ is drawn as in Figure 4.5. 
\f
For $0\leq\theta\leq\pi$, $F_\theta\cap(L\cup h')$ is drawn as in Figure 4.2. 
\f
(3) $B^4\cap h'$ satisfies the following conditions. 
\f
For $\pi<\theta<2\pi$, $F_\theta\cap h'$ is empty. 
\f
For $0\leq\theta\leq\pi$, $F_\theta\cap h'$ is drawn as in Figure 4.3. 
\f
(4)
Note that $h$ and $h'$ are dual handles each other. 
\f
(5)
$B^4\cap K_3$ satisfies the following conditions. 
(We summarize the conditions in Table 1.)
\f
For $\pi<\theta<2\pi$, $F_\theta\cap K_3$ is drawn as in Figure 4.5. 
\f
For $\theta=0,\pi$, $F_\theta\cap K_3$ is drawn as in Figure 4.2.  
\f
For $0<\theta<\pi$, $F_\theta\cap K_3$ is drawn as in Figure 4.1.

By  the construction of  this knot $K_3$ and the construction of  $Q(K)$ 
in \S4, $K_3$ is identical to $Q(K)$. 
This completes the proof of Lemma 4.1, Proposition 3.1, Theorem 1.1.

\newpage
\unitlength 0.1in
\begin{picture}(53.80,48.20)(2.00,-47.80)
%
\special{pn 8}%
\special{pa 1000 790}%
\special{pa 1000 4780}%
\special{fp}%
%
\special{pn 8}%
\special{pa 2010 30}%
\special{pa 1000 810}%
\special{fp}%
%
\special{pn 8}%
\special{pa 2010 3990}%
\special{pa 1000 4770}%
\special{fp}%
%
\special{pn 8}%
\special{pa 2010 3990}%
\special{pa 5530 3990}%
\special{fp}%
%
\special{pn 8}%
\special{pa 1600 1080}%
\special{pa 3700 1080}%
\special{fp}%
%
\special{pn 8}%
\special{pa 1610 800}%
\special{pa 4010 800}%
\special{fp}%
\put(7.9000,-47.1000){\makebox(0,0)[lb]{y}}%
\put(20.9000,-1.3000){\makebox(0,0)[lb]{x}}%
\put(19.0000,-3.6000){\makebox(0,0)[lb]{1}}%
\put(10.7000,-45.9000){\makebox(0,0)[lb]{1}}%
\put(54.0000,-42.1000){\makebox(0,0)[lb]{1}}%
\put(20.1000,-41.5000){\makebox(0,0)[lb]{o}}%
%
\special{pn 8}%
\special{pa 3710 1280}%
\special{pa 2690 1280}%
\special{fp}%
%
\special{pn 8}%
\special{pa 2700 1290}%
\special{pa 2700 2140}%
\special{fp}%
%
\special{pn 8}%
\special{pa 4820 2130}%
\special{pa 4820 2490}%
\special{fp}%
%
\special{pn 8}%
\special{pa 4820 2510}%
\special{pa 2700 2500}%
\special{fp}%
%
\special{pn 8}%
\special{pa 2710 2500}%
\special{pa 2710 3260}%
\special{fp}%
%
\special{pn 8}%
\special{pa 4820 3290}%
\special{pa 4820 3530}%
\special{fp}%
%
\special{pn 8}%
\special{pa 4820 3510}%
\special{pa 4020 3530}%
\special{fp}%
%
\special{pn 8}%
\special{pa 3710 1090}%
\special{pa 3710 1280}%
\special{fp}%
%
\special{pn 8}%
\special{pa 4050 3530}%
\special{pa 4050 2910}%
\special{fp}%
%
\special{pn 8}%
\special{pa 4060 2680}%
\special{pa 4050 2580}%
\special{fp}%
%
\special{pn 8}%
\special{pa 4050 2380}%
\special{pa 4040 800}%
\special{fp}%
%
\special{pn 8}%
\special{pa 3980 810}%
\special{pa 4040 800}%
\special{fp}%
%
\special{pn 8}%
\special{pa 2700 2140}%
\special{pa 3570 2140}%
\special{fp}%
%
\special{pn 8}%
\special{pa 3800 2130}%
\special{pa 3960 2120}%
\special{fp}%
%
\special{pn 8}%
\special{pa 4190 2120}%
\special{pa 4830 2140}%
\special{fp}%
%
\special{pn 8}%
\special{pa 380 650}%
\special{pa 390 650}%
\special{fp}%
%
\special{pn 8}%
\special{pa 200 1200}%
\special{pa 200 1210}%
\special{fp}%
%
\special{pn 8}%
\special{pa 2710 3260}%
\special{pa 3600 3270}%
\special{fp}%
%
\special{pn 8}%
\special{pa 3800 3290}%
\special{pa 3970 3300}%
\special{fp}%
%
\special{pn 8}%
\special{pa 4170 3280}%
\special{pa 4820 3310}%
\special{fp}%
\put(41.8000,-7.9000){\makebox(0,0)[lb]{l}}%
\put(42.5000,-8.5000){\makebox(0,0)[lb]{1}}%
%
\special{pn 8}%
\special{pa 3580 2140}%
\special{pa 3800 2130}%
\special{fp}%
%
\special{pn 8}%
\special{pa 3610 3270}%
\special{pa 3790 3290}%
\special{fp}%
%
\special{pn 8}%
\special{pa 4060 2620}%
\special{pa 4050 2990}%
\special{fp}%
\put(25.7000,-47.6000){\makebox(0,0)[lb]{Figure 5.1}}%
\put(55.8000,-39.2000){\makebox(0,0)[lb]{r}}%
%
\special{pn 8}%
\special{pa 2010 60}%
\special{pa 2010 30}%
\special{fp}%
%
\special{pn 8}%
\special{pa 2010 40}%
\special{pa 2010 710}%
\special{fp}%
%
\special{pn 8}%
\special{pa 2010 880}%
\special{pa 2010 1020}%
\special{fp}%
%
\special{pn 8}%
\special{pa 2010 1180}%
\special{pa 2010 3980}%
\special{fp}%
\end{picture}%

\newpage
\unitlength 0.1in
\begin{picture}(53.80,48.20)(2.00,-47.80)
%
\special{pn 8}%
\special{pa 2010 10}%
\special{pa 2010 4000}%
\special{fp}%
%
\special{pn 8}%
\special{pa 1000 790}%
\special{pa 1000 4780}%
\special{fp}%
%
\special{pn 8}%
\special{pa 2010 30}%
\special{pa 1000 810}%
\special{fp}%
%
\special{pn 8}%
\special{pa 2010 3990}%
\special{pa 1000 4770}%
\special{fp}%
\put(7.9000,-47.1000){\makebox(0,0)[lb]{y}}%
\put(20.9000,-1.3000){\makebox(0,0)[lb]{x}}%
\put(19.0000,-3.6000){\makebox(0,0)[lb]{1}}%
\put(10.7000,-45.9000){\makebox(0,0)[lb]{1}}%
\put(54.0000,-42.1000){\makebox(0,0)[lb]{1}}%
\put(20.1000,-41.5000){\makebox(0,0)[lb]{o}}%
%
\special{pn 8}%
\special{pa 3710 1530}%
\special{pa 2910 1530}%
\special{fp}%
%
\special{pn 8}%
\special{pa 2920 1530}%
\special{pa 2920 1890}%
\special{fp}%
%
\special{pn 8}%
\special{pa 5110 1910}%
\special{pa 5110 2750}%
\special{fp}%
%
\special{pn 8}%
\special{pa 5110 2760}%
\special{pa 2920 2760}%
\special{fp}%
%
\special{pn 8}%
\special{pa 2940 2760}%
\special{pa 2940 3050}%
\special{fp}%
%
\special{pn 8}%
\special{pa 5110 3050}%
\special{pa 5110 3750}%
\special{fp}%
%
\special{pn 8}%
\special{pa 5110 3780}%
\special{pa 4020 3790}%
\special{fp}%
%
\special{pn 8}%
\special{pa 3710 1540}%
\special{pa 3710 2390}%
\special{fp}%
%
\special{pn 8}%
\special{pa 3710 4360}%
\special{pa 3710 2900}%
\special{fp}%
%
\special{pn 8}%
\special{pa 3710 2570}%
\special{pa 3710 2680}%
\special{fp}%
%
\special{pn 8}%
\special{pa 4040 3780}%
\special{pa 4040 4360}%
\special{fp}%
%
\special{pn 8}%
\special{pa 2920 1880}%
\special{pa 3570 1880}%
\special{fp}%
%
\special{pn 8}%
\special{pa 3800 1880}%
\special{pa 3970 1880}%
\special{fp}%
%
\special{pn 8}%
\special{pa 4190 1890}%
\special{pa 5120 1910}%
\special{fp}%
%
\special{pn 8}%
\special{pa 2930 3040}%
\special{pa 3600 3040}%
\special{fp}%
%
\special{pn 8}%
\special{pa 3800 3040}%
\special{pa 3980 3040}%
\special{fp}%
%
\special{pn 8}%
\special{pa 4170 3030}%
\special{pa 5110 3060}%
\special{fp}%
%
\special{pn 8}%
\special{pa 380 650}%
\special{pa 390 650}%
\special{fp}%
%
\special{pn 8}%
\special{pa 200 1200}%
\special{pa 200 1210}%
\special{fp}%
\put(52.0000,-35.9000){\makebox(0,0)[lb]{l}}%
\put(52.6000,-36.7000){\makebox(0,0)[lb]{2}}%
%
\special{pn 8}%
\special{pa 5110 3790}%
\special{pa 5110 3720}%
\special{fp}%
%
\special{pn 8}%
\special{pa 4000 3040}%
\special{pa 4240 3030}%
\special{fp}%
%
\special{pn 8}%
\special{pa 4050 3050}%
\special{pa 3950 3050}%
\special{fp}%
%
\special{pn 8}%
\special{pa 3960 1890}%
\special{pa 4210 1880}%
\special{fp}%
%
\special{pn 8}%
\special{pa 3710 2650}%
\special{pa 3710 2340}%
\special{fp}%
\put(25.0000,-47.9000){\makebox(0,0)[lb]{Fgure 5.2}}%
\put(55.8000,-39.2000){\makebox(0,0)[lb]{r}}%
%
\special{pn 8}%
\special{pa 2010 4000}%
\special{pa 3570 4000}%
\special{fp}%
%
\special{pn 8}%
\special{pa 3780 4000}%
\special{pa 3990 4000}%
\special{fp}%
%
\special{pn 8}%
\special{pa 4110 4000}%
\special{pa 5520 4000}%
\special{fp}%
\end{picture}%

\newpage

\section{  Related problems }

\noindent
{\bf Problem 6.1.}
Let $L=(K_1, K_2)$ be a 2-link. Then do we have $\mu(L)=\mu(K_1)+\mu(K_2)$? 
\vskip3mm

See 
\cite{R}  
\cite{O3} 
\cite{O4} 
for the $\mu$ invariant of 2-links and related topics.

\vskip3mm
\noindent
{\bf Problem 6.2.}
Is there a 2-link which is not an SHB link? 

\vskip3mm
See 
\cite{CO} 
\cite{LO}
for SHB links.

The author proved in 
\cite{O3}:   
if  $L=(K_1, K_2)$ is an SHB link, 
then the answer to Problem 6.1 is affirmative.

\vskip3mm
\noindent
{\bf Problem 6.3.}
Let $K_1, K_2, K_3$ be arbitrary 2-knots. 
Is there a 2-component 2-link $L=(L_1,L_2)$ such that 
$L_1$ (resp. $L_2$) is equivalent to $K_2$ (resp. $K_2$) 
and that a band-sum of $L$ is $K_3$? 
\vskip3mm

If the answer to Problem 6.1 is affirmative, 
then the answer to Problem 6.3 is negative. 

In 
\cite{O2}  
the author gives 
the negative answer to the $n$-dimensional knot version of Problem 6.3  
and proves  the $n$-dimensional version of Theorem 1.1.  
The announcement of them is in 
\cite{O1}.






{\bf Acknowledgements.  } 
The author discussed Problem 6.1 with Prof. Kent Orr. 
Prof. Kent Orr, Prof. J. Levine and Prof. T. Cochran 
told the author that Problem 6.2 is open. 
The author would like to thank  Prof. Kent Orr, 
Prof. J. Levine and Prof. T. Cochran for their kindness.
The author would like to thank the referee(s) and the editor(s)
for reading the manuscript. 

\end{document}